\documentclass[12pt,thmsa,oneside]{article}
\usepackage{amsfonts}
\usepackage{sw20lart}

\input tcilatex
\QQQ{Language}{
American English
}

\begin{document}

\pagestyle{myheadings} \markright{} \TeXButton{title}
{\title{
\vspace{-1.5in}
\begin{flushright}
{\small UAHEP 9702}
\end{flushright}
\vspace{-48pt}
\begin{flushright}
{\small February 1997}
\end{flushright}
\vspace{24pt}
Towards constructing one-particle representations of the deformed Poincar\'e algebra
}}

\TeXButton{author}
{\author{I. Yakushin
\thanks{  \hspace{.1in} Department of Physics and Astronomy, 
University of Alabama, Tuscaloosa, USA}}}

\TeXButton{date}{\date{}}

\TeXButton{maketitle}{\maketitle}

\TeXButton{babstract}{\begin{abstract}}We give a method for obtaining states
of massive particle representations of the two-parameter deformation of the
Poincar\'e algebra proposed in \cite{sy96a}, \cite{sy96}, \cite{ssy95}. We
discuss four procedures to generate eigenstates of a complete set of
commuting operators starting from the rest state. One result of this work is
the fact that upon deforming to the quantum Poincar\'e algebra the rest
state is split into an infinite number of states. Another result is that the
energy spectrum of these states is discrete. Some curious residual
degeneracy remains: there are states constructed by applying different
operators to the rest state which nevertheless are indistinguishable by
eigenvalues of all the observables in the algebra.

\TeXButton{eabstract}{\end{abstract}}

\section{Introduction\label{int4}\ }

Various deformations of the Poincar\'e algebra have been given in the
literature (\cite{oswz92}, \cite{PSW93}, \cite{swz91}, \cite{sm1194}, \cite
{sm895}, \cite{nl91}, \cite{rt94}, \cite{lr94}, \cite{bhos95}). These
deformations are speculated to play a role at Planck length scale. In this
article we give procedures for constructing some states of non-zero mass
representations of the two-parameter deformation of the Poincar\'e algebra
which we proposed in \cite{sy96a}, \cite{sy96}, \cite{ssy95}.

Let us first briefly review the results of \cite{sy96a}, \cite{sy96}, \cite
{ssy95}. In \cite{sy96}, \cite{ssy95} we discussed a Lie-Poisson deformation
of the Poincar\'e algebra. For such an algebra the Lorentz group $SL(2,C)$
does not act canonically, but instead as a Lie-Poisson group. This algebra
depends on one deformation parameter $\lambda .$ In \cite{sy96} we showed
how to quantize this algebra which, of course, introduces another parameter $%
\hbar $ that also can be regarded as a deformation parameter. In the limit $%
\lambda \rightarrow 0$ we recover the standard Poincar\'e algebra from the
quantum algebra, while the limit $\hbar \rightarrow 0$ gives us the
Lie-Poisson deformation introduced in \cite{ssy95}. In \cite{sy96a} we
investigated properties of this two-parameter deformation of the Poincar\'e
algebra. We found that it is covariant under $SL_q(2,C),$ determined its
Casimirs and obtained the complete set of commuting operators. We also found
a subalgebra which is a curious deformation of $su(2).$

Next we summarize the results of \cite{sy96a} which are necessary to
construct a representation. The algebra can be compactly written as follows: 
\begin{eqnarray}
\stackunder{12}{R}\stackunder{1}{P}\stackunder{12}{R}^{-1}\stackunder{2}{P}
&=&\stackunder{2}{P}\stackunder{21}{R}^{-1}\stackunder{1}{P}\stackunder{21}{R%
},  \nonumber  \label{pp4} \\
\stackunder{21}{R}^{-1}\stackunder{1}{\Gamma }\stackunder{21}{R}\stackunder{2%
}{\Gamma } &=&\stackunder{2}{\Gamma }\stackunder{12}{R}\stackunder{1}{\Gamma 
}\stackunder{12}{R}^{-1},  \nonumber  \label{gg4} \\
\stackunder{12}{R}\stackunder{1}{\Gamma }\stackunder{12}{R}^{-1}\stackunder{2%
}{\overline{\Gamma }} &=&\stackunder{2}{\overline{\Gamma }}\stackunder{12}{R}%
\stackunder{1}{\Gamma }\stackunder{12}{R}^{-1},  \label{pggb4} \\
\stackunder{21}{R}^{-1}\stackunder{1}{P}\stackunder{21}{R}\stackunder{2}{%
\Gamma } &=&\stackunder{2}{\Gamma }\stackunder{21}{R}^{-1}\stackunder{1}{P}%
\stackunder{12}{R}^{-1},  \nonumber  \label{pg4}
\end{eqnarray}
We use tensor product notation, labels 1 and 2 denote different vector
spaces, $\stackunder{1}{P}\equiv P\otimes \Bbb{I},$ $\stackunder{2}{P}\equiv 
\Bbb{I}\otimes P$ , etc., $\Bbb{I}$ is 2$\times $2 unit matrix. $P,$ $\Gamma 
$ and $\overline{\Gamma }$ are $2\times 2$ matrices. $P$ contains the four
momentum operators, while $\Gamma $ and $\overline{\Gamma }$ contain the
angular momentum operators. Relations (\ref{pggb4}) are written out in
components in the Appendix. The 2$\times $2 matrices satisfy the following
properties:

\begin{enumerate}
\item[i)]  $P$ is assumed to be hermitian on the states of the representation%
\footnote{%
From now on we might use hermitian conjugation without mentioning the
representation but, of course, without representation it is not defined, we
simply do not want to write each time ''on the states of our
representation''.} we are constructing and can be expressed via components $%
P_\alpha $ of the momentum 4-vector as follows: 
\begin{eqnarray}
P &=&\left( 
\begin{array}{cc}
P_{11} & P_{12} \\ 
P_{21} & P_{22}
\end{array}
\right) =-\Bbb{I}P_0+\sigma _kP_k  \label{P4} \\
&=&\left( 
\begin{array}{cc}
-P_0+P_3 & P_1-i~P_2 \\ 
P_1+i~P_2 & -P_0-P_3
\end{array}
\right) =P^{\dagger },  \nonumber \\
P_\alpha &=&P_\alpha ^{\dagger },\quad \alpha =0,1,2,3.
\end{eqnarray}

\item[ii)]  $\Gamma $ and $\overline{\Gamma }$ satisfy a deformed
unimodularity condition: 
\begin{eqnarray}
det_{\frac 1q}(\Gamma ^T)\equiv \Gamma _{11}\Gamma _{22}-q^2\Gamma
_{21}\Gamma _{12}=1,  \nonumber  \label{dg4} \\
det_q(\overline{\Gamma })\equiv \overline{\Gamma }_{11}\overline{\Gamma }%
_{22}-\frac 1{q^2}\overline{\Gamma }_{12}\overline{\Gamma }_{21}=1.
\label{dggb4}
\end{eqnarray}

\item[iii)]  $\Gamma $ and $\overline{\Gamma }$ are presumed to satsify
relations 
\begin{eqnarray}
\overline{\Gamma }^{-1} &=&\Gamma ^{\dagger },  \label{ggbdag4} \\
\Gamma ^{-1} &=&\overline{\Gamma }^{\dagger }  \nonumber
\end{eqnarray}
and can be expressed in terms of angular momentum tensor $J$ as follows: 
\[
\Gamma =e^{i\lambda J},\quad \overline{\Gamma }=e^{i\lambda J^{\dagger }}, 
\]
where $\lambda $ is a real parameter and 
\begin{eqnarray}
J &=&\sigma _k\left( J_{k0}-\frac i2\epsilon _{kmn}J_{mn}\right) =
\label{J4} \\
&=&\left( 
\begin{array}{cc}
-iJ_{12}+J_{30} & -iJ_{23}-iJ_{20}-J_{31}+J_{10} \\ 
-iJ_{23}+iJ_{20}+J_{31}+J_{10} & iJ_{12}-J_{30}
\end{array}
\right) ,  \nonumber \\
J_{\alpha \beta } &=&J_{\alpha \beta }^{\dagger },\quad \alpha ,\beta
=0,1,2,3.  \nonumber
\end{eqnarray}
\end{enumerate}

The $R-$ matrix is given by 
\begin{equation}
\stackunder{12}{R}=q^{-1/2}\left( 
\begin{array}{cccc}
q & 0 & 0 & 0 \\ 
0 & 1 & 0 & 0 \\ 
0 & q-q^{-1} & 1 & 0 \\ 
0 & 0 & 0 & q
\end{array}
\right)  \label{R4}
\end{equation}
and satisfies the quantum Yang-Baxter equation. Here $q=e^{\hbar \lambda }\ $
and $\hbar $ and $\lambda $ can be regarded as deformation parameters. In 
\cite{sy96a} we showed that in the limit $\lambda \rightarrow 0$ we obtain
the usual Poincar\'e algebra, while the limit $\hbar \rightarrow 0$ gives us
Lie-Poisson deformation of the Poincar\'e algebra which we discussed in \cite
{sy96} and \cite{ssy95}. As was shown in \cite{sy96a}, the algebra is
preserved under the action of $SL_q(2,C):$%
\begin{eqnarray}
P &\rightarrow &P^{\prime }=\overline{T}PT^{-1},  \label{sym4} \\
\Gamma &\rightarrow &\Gamma ^{\prime }=T\Gamma T^{-1},  \nonumber \\
\overline{\Gamma } &\rightarrow &\overline{\Gamma }^{\prime }=\overline{T}~%
\overline{\Gamma }~\overline{T}^{-1}.  \nonumber
\end{eqnarray}
where $T\in SL_q(2,C),\quad T^{\dagger }=\overline{T}^{-1}.$ By definition $%
T $ satisfies the commutational relations 
\begin{eqnarray}
\stackunder{12}{R}\stackunder{1}{T}\stackunder{2}{T} &=&\stackunder{2}{T}%
\stackunder{1}{T}\stackunder{12}{R},  \nonumber  \label{att4} \\
\stackunder{12}{R}\stackunder{1}{T}\stackunder{2}{\overline{T}} &=&%
\stackunder{2}{\overline{T}}\stackunder{1}{T}\stackunder{12}{R},  \label{tt4}
\\
\stackunder{12}{R}~\stackunder{1}{\overline{T}}\stackunder{2}{~\overline{T}}
&=&\stackunder{2}{\overline{T}}\stackunder{1}{~\overline{T}}~\stackunder{12}{%
R},  \nonumber
\end{eqnarray}
as well as the deformed unimodularity condition 
\begin{equation}
det_{\frac 1{\sqrt{q}}}\left( T\right) \equiv T_{11}T_{22}-qT_{12}T_{21}=1
\end{equation}
Our algebra appears to be distinct from systems discussed previously, e.g.
in \cite{oswz92}, \cite{PSW93}, \cite{swz91}, \cite{sm1194}, \cite{sm895}, 
\cite{nl91}, \cite{rt94}, \cite{lr94}, \cite{bhos95} and can be expressed
compactly.

The following matrix was shown \cite{sy96a} to be an analogue of
Pauli-Lubanski vector: 
\begin{equation}
W\equiv -\Bbb{I}W_0+\sigma _kW_k=a\left( \beta \overline{\Gamma }%
^{-1}P\Gamma -P\right) ,  \label{w4}
\end{equation}
where $W^{\dagger }=W,$ and $a$ and $\beta $ are c-numbers with the
following limits: $a\stackunder{\hbar \rightarrow 0}{\rightarrow }\frac
1{2\lambda },$ $\beta \stackunder{\hbar \rightarrow 0}{\rightarrow }1,$ $a%
\stackunder{\lambda \rightarrow 0}{\rightarrow }\frac 1{2\lambda },$ $\beta 
\stackunder{\lambda \rightarrow 0}{\rightarrow }1.$ We shall see in Sec. \ref
{eigen4} that $\beta =q^3$ .

The algebra has two Casimir operators which are also invariant under $%
SL_q(2,C):$%
\begin{eqnarray}
\mathcal{C}_1 &=&(P,P)_q~,  \label{Casimirs4} \\
~\mathcal{C}_2 &=&(W,W)_q,  \nonumber
\end{eqnarray}
where the deformed scalar product is defined as follows: 
\begin{equation}
(A,B)_q=-\frac 1{q^2+1}Tr_q\left( A\widetilde{B}\right) ,  \label{()4}
\end{equation}
\[
Tr_q(A)\equiv A_{11}+q^2A_{22}, 
\]
The twiddle denotes an adjugate. It is defined as follows: if $B$ is matrix,
then its adjugate $\widetilde{B}$ has the properties 
\[
B\widetilde{B}=\widetilde{B}B\sim \Bbb{I} 
\]
and in the limit $q\rightarrow 1$ it goes to $B^{-1}\det \left( B\right) .$ $%
\widetilde{B}$ is defined uniquely up to a constant which goes to $1$ in the
limit $q\rightarrow 1.$ For the matrices $P$ and $W$ we find 
\begin{equation}
\widetilde{P}=\left( 
\begin{array}{cc}
P_{22} & -\frac 1{q^2}P_{12} \\ 
-\frac 1{q^2}P_{21} & \frac 1{q^2}\left( P_{11}+(q^2-1)P_{22}\right)
\end{array}
\right) .  \label{Pad4}
\end{equation}
\begin{equation}
\widetilde{W}=\left( 
\begin{array}{cc}
W_{22} & -\frac{W_{12}}{q^2} \\ 
-\frac{W_{21}}{q^2} & \frac{W_{11}+(q^2-1)W_{22}}{q^2}
\end{array}
\right) -a(q^2-1)\widetilde{P}  \label{Wad4}
\end{equation}
When $\lambda \rightarrow 0,$ we recover the usual form for the
Pauli-Lubanski vector and Casimir operators: 
\[
W_\beta \stackunder{\lambda \rightarrow 0}{\rightarrow }-\frac 12\epsilon
_\beta ^{~\mu \nu \rho }J_{\mu \nu }P_\rho ,\quad \mathcal{C}_1\stackunder{%
\lambda \rightarrow 0}{\rightarrow }P_\mu P^\mu ,\quad \mathcal{C}_2%
\stackunder{\lambda \rightarrow 0}{\rightarrow }W_\mu W^\mu . 
\]

From the relations $\left( \text{\ref{P4}}\right) $ and $\left( \text{\ref
{ggbdag4}}\right) $ we are assuming that matrices $P,$ $\Gamma +\overline{%
\Gamma }^{-1},$ $\overline{\Gamma }+\Gamma ^{-1}$ are hermitian\footnote{%
where 
\begin{eqnarray*}
\overline{\Gamma }^{-1} &=&\left( 
\begin{array}{cc}
\overline{\Gamma }_{22} & -\frac 1{q^2}\overline{\Gamma }_{12} \\ 
-\frac 1{q^2}\overline{\Gamma }_{21} & \frac 1{q^2}\left( \overline{\Gamma }%
_{11}+(q^2-1)\overline{\Gamma }_{22}\right)
\end{array}
\right) ,  \label{gbinv3} \\
\Gamma ^{-1} &=&\left( 
\begin{array}{cc}
q^2\Gamma _{22}-(q^2-1)\Gamma _{11} & -q^2\Gamma _{12} \\ 
-q^2\Gamma _{21} & \Gamma _{11}
\end{array}
\right) .  \label{ginv3}
\end{eqnarray*}
} on the states of our representations. Each hermitian matrix $A$ gives 4
hermitian operators $A_{11},$ $A_{22},$ $A_{12}+A_{21},$ $\frac{A_{12}-A_{21}%
}i.$ Therefore we have 12 hermitian operators and two constraints $\left( 
\text{\ref{dggb4}}\right) .$ Hence we can construct 10 independent hermitian
linear combinations of the generators of our algebra. This means that by
construction our representation is unitary. We will find no contradictions
to this assumption later, in particular all the eigenvalues of the
physically meaningful operators are real.

In section \ref{su24}, \ref{lo4}, \ref{eu4} we examine three different
subalgebras of the full algebra generated by $\Gamma ,$ $\overline{\Gamma }$
and $P.$ The first is a deformation of the 3-angular momentum algebra. Its
representations can be easily constructed and are isomorphic to the
representation of $su(2)$ . In sections \ref{lo4} and \ref{eu4} we consider
the analogues of Lorentz and Euclidean subalgebras. These subalgebras are
more involved and we only give the Casimir operators for them. In section 
\ref{set4} we write out the complete set of commuting operators. Finally, in
section \ref{eigen4} we find four classes of eigenstates of all the
operators in the complete set of commuting operators. It is shown that in
the limit $\lambda \rightarrow 0$ all the eigenvalues correspond to those of
the rest state of the non-deformed Poincar\'e algebra. Concluding remarks
are made in section \ref{concl4}.

\section{Representations of the deformed su(2) subalgebra\label{su24}}

The central role in our construction is played by the universal enveloping
algebra generated by the hermitian matrix 
\begin{equation}
\Omega \equiv \Gamma \overline{\Gamma }^{-1}=\Gamma \Gamma ^{\dagger }
\label{omega4}
\end{equation}
which has very surprising property of commuting in the same way with all the
matrices defined above: 
\begin{equation}
\stackunder{1}{Z}\stackunder{21}{R}\stackunder{2}{\Omega }\stackunder{12}{R}=%
\stackunder{21}{R}\stackunder{2}{\Omega }\stackunder{12}{R}\stackunder{1}{Z},
\label{sigmacr4}
\end{equation}
where $Z=P,$ $\Gamma ,$ $\overline{\Gamma },$ $W$ or $\Omega .$ Not all four
components of $\Omega $ are independent. From the properties of $\Gamma $
and $\overline{\Gamma }$ it follows that 
\begin{equation}
det_{\frac 1q}(\Omega ^T)\equiv \Omega _{11}\Omega _{22}-q^2\Omega
_{21}\Omega _{12}=1.  \label{sigmaconstr4}
\end{equation}

We shall show that the representations of the algebra of $\Omega _{ij}$ are
in one-to one correspondence with the representations of the $su(2)$
algebra. Furthermore, from the $\lambda \rightarrow 0$ limit 
\begin{equation}
\frac{\Omega -\Bbb{I}}{2\lambda }\stackunder{\lambda \rightarrow 0}{%
\rightarrow }\left( 
\begin{array}{ll}
J_{12} & J_{23}-iJ_{31} \\ 
J_{23}+iJ_{31} & -J_{12}
\end{array}
\right)  \label{omegalim4}
\end{equation}
we can see that the $\Omega -$ algebra is a deformation of the $su(2)$ -
algebra with $\Omega _{11}$ being analogous to the projection of 3-vector of
angular momentum onto the third axis. The Casimir of the $\Omega -$ algebra
is 
\begin{equation}
Tr_q\left( \Omega \right) \stackunder{\lambda \rightarrow 0}{\rightarrow }%
2+2\hbar \lambda +\left( \hbar ^2-4i\hbar J_{30}+4\left(
J_{12}^2+J_{23}^2+J_{31}^2\right) \right) \lambda ^2+o\left( \lambda
^2\right)  \label{k54}
\end{equation}
which corresponds to the square of 3-vector of angular momentum. [We note
that although $Tr_q\left( \Omega \right) $ commutes with $\Omega _{ij},$ it
does not in general commute with $P_{ij}$, $\Gamma _{ij}$and $\overline{%
\Gamma }_{ij}.$ ] Below we shall look for states which are diagonal in $%
\Omega _{11}$ and $Tr_q\left( \Omega \right) .$

Let us rewrite the $\Omega -$ algebra (\ref{sigmacr4}) [with $Z=\Omega $ ]
in components: 
\begin{eqnarray}
\Omega _{12}\Omega _{11} &=&q^2\Omega _{11}\Omega _{12},\quad  \nonumber
\label{sigmacomp4a} \\
\Omega _{21}\Omega _{11} &=&\frac 1{q^2}\Omega _{11}\Omega _{21},\quad 
\nonumber \\
\Omega _{22}\Omega _{11} &=&\Omega _{11}\Omega _{22},  \nonumber \\
\Omega _{22}\Omega _{12} &=&\Omega _{12}\Omega _{22}+\frac{q^2-1}{q^4}\Omega
_{12}\Omega _{11},\quad  \nonumber  \label{sigmacomp4b} \\
\Omega _{21}\Omega _{12} &=&\Omega _{12}\Omega _{21}-\frac{q^2-1}{q^2}\left(
\Omega _{22}-\Omega _{11}\right) \Omega _{11},  \nonumber \\
\Omega _{22}\Omega _{21} &=&\Omega _{21}\Omega _{22}-\frac{q^2-1}{q^2}\Omega
_{21}\Omega _{11},\quad  \label{sigmacomp4} \\
\Omega _{11}\Omega _{22}-q^2\Omega _{21}\Omega _{12} &=&1.  \nonumber
\end{eqnarray}
In addition the condition that $\Omega $ is hermitian means that 
\[
\Omega _{11}=\Omega _{11}^{\dagger },\quad \Omega _{22}=\Omega
_{22}^{\dagger },\quad \Omega _{12}=\Omega _{21}^{\dagger }. 
\]
One can check that $Tr_q(\Omega )$ and $\Omega _{11}$ form a complete set of
commuting operators in the $\Omega -$ algebra. From the first two
commutational relations (\ref{sigmacomp4}) it is obvious that: 
\begin{eqnarray*}
Tr_q(\Omega )\left| k,m\right\rangle &=&k\left| k,m\right\rangle , \\
\Omega _{11}\left| k,m\right\rangle &=&\rho q^{2m}\left| k,m\right\rangle ,
\\
\Omega _{12}\left| k,m\right\rangle &=&A_{k,m}\left| k,m-1\right\rangle , \\
\Omega _{21}\left| k,m\right\rangle &=&B_{k,m}\left| k,m+1\right\rangle .
\end{eqnarray*}
For simplicity we shall assume $\rho =1.$ From hermiticity of $\Omega :$%
\[
A_{k,m}=B_{k,m-1}^{*} 
\]
From $det_{\frac 1q}(\Omega ^T)=1:$%
\begin{eqnarray*}
\left| A_{k,m}\right| ^2 &=&\frac{q^{2m}k-q^2-q^{4m}}{q^4}, \\
\left| B_{k,m}\right| ^2 &=&\frac{q^{2m}k-1-q^{4m+2}}{q^2}.
\end{eqnarray*}
In order $\left| B_{k,m}\right| ^2\geq 0,$ one must require 
\[
k\geq q\left( q^{2m+1}+q^{-2m-1}\right) . 
\]
Since the right hand side monotonously increases for positive $m,$ there
exists a maximum value of $m$ which we denote by $j$ such that $B_{k,j}=0$
and $k=q\left( q^{2j+1}+q^{-2j-1}\right) .$ Similar logic leads to the
conclusion that there exists a minimum value of $m$ and finally 
\[
-j\leq m\leq j. 
\]
It is therefore more convenient to use $j$ and $m$ to label the
representation. We shall also use below the following function: 
\begin{equation}
k_j\equiv q\left( q^{2j+1}+q^{-\left( 2j+1\right) }\right) .  \label{ki4}
\end{equation}
Let us rewrite our representation in this language: 
\begin{eqnarray}
Tr_q(\Omega )\left| j,m\right\rangle &=&k_j\left| j,m\right\rangle , 
\nonumber  \label{omegarep4a} \\
\Omega _{11}\left| j,m\right\rangle &=&q^{2m}\left| j,m\right\rangle , 
\nonumber \\
\Omega _{12}\left| j,m\right\rangle &=&A_{j,m}\left| j,m-1\right\rangle ,
\label{omegarep4} \\
\Omega _{21}\left| j,m\right\rangle &=&B_{j,m}\left| j,m+1\right\rangle , 
\nonumber \\
\left| B_{j,m}\right| ^2 &=&q^{2\left( m-1\right) }\left( k_j-k_m\right)
,\quad  \nonumber \\
A_{j,m} &=&B_{j,m-1}^{*}.  \nonumber
\end{eqnarray}
To see how $k_j$ is analogous to $j(j+1),$ consider the limit: 
\[
a^2\left( k_j-k_0\right) \stackunder{\lambda \rightarrow 0}{\rightarrow }%
\hbar ^2j(j+1), 
\]
where $a$ is the same as in $\left( \text{\ref{w4}}\right) .$

Just as with the $su(2)$ representations, $j$ can be either an integer or
half-integer. Since $m$ changes by 1 between $j$ and $-j$ there should exist
an integer $n$ such that $j-n=-j,$ that is $j=\frac n2.$

For $j=1/2$, for example, we can write out the representation of $\Omega
^{\prime }s$ in terms of $2\times 2$ matrices: 
\begin{eqnarray*}
\Omega _{11} &=&\left( 
\begin{array}{ll}
q & 0 \\ 
0 & \frac 1q
\end{array}
\right) ,~~~~~~~\Omega _{12}=\left( 
\begin{array}{ll}
0 & 0 \\ 
\frac{q^2-1}{q^2} & 0
\end{array}
\right) ,~ \\
\Omega _{21} &=&\left( 
\begin{array}{ll}
0 & \frac{q^2-1}{q^2} \\ 
0 & 0
\end{array}
\right) ,~~~\Omega _{22}=\left( 
\begin{array}{ll}
\frac{q^2\left( q^2-1\right) +1}{q^3} & 0 \\ 
0 & q
\end{array}
\right)
\end{eqnarray*}
assuming that $\left| \frac 12,\frac 12\right\rangle =\left( 
\begin{array}{l}
1 \\ 
0
\end{array}
\right) ,$ $\left| \frac 12,-\frac 12\right\rangle =\left( 
\begin{array}{l}
0 \\ 
1
\end{array}
\right) $ as usual. All commutational relations and constraints are
satisfied for this representation.

\section{Deformed Lorentz subalgebra\label{lo4}}

The deformed Lorentz subalgebra is generated by $\Gamma $ and $\overline{%
\Gamma }$ and contains the $\Omega -$ subalgebra. In the non-deformed case
there are two Casimirs: $\left( J_{12}^2+J_{23}^2+J_{31}^2\right) -\left(
J_{10}^2+J_{20}^2+J_{30}^2\right) $ and $%
J_{12}J_{30}+J_{23}J_{10}+J_{31}J_{20}.$ In the complete set of commuting
operators one can include $J_{12}$ and $J_{30}$ with the Casimirs. In our
algebra $Tr_q\left( \Gamma +\overline{\Gamma }\right) $ and $Tr_q\left(
\Gamma -\overline{\Gamma }\right) $ are Casimirs of the deformed Lorentz
subalgebra. From their limits 
\begin{eqnarray*}
&&\ \ \ Tr_q\left( \Gamma +\overline{\Gamma }\right) \stackunder{\lambda
\rightarrow 0}{\rightarrow }  \label{trg+gb} \\
&&\ \ \ 4\left( 1+\hbar \lambda +\frac 12\left( \hbar ^2-2i\hbar
J_{30}+\left( J_{12}^2+J_{23}^2+J_{31}^2\right) -\left(
J_{10}^2+J_{20}^2+J_{30}^2\right) \right) \lambda ^2\right)  \label{trg-gb}
\\
&&\ +O\left( \lambda ^3\right) , \\
&&\ Tr_q\left( \Gamma -\overline{\Gamma }\right) \stackunder{\lambda
\rightarrow 0}{\rightarrow }4i\left( i\hbar
J_{12}+J_{12}J_{30}+J_{23}J_{10}+J_{31}J_{20}\right) \lambda ^2+O\left(
\lambda ^3\right)
\end{eqnarray*}
we conclude that they correspond to $\left(
J_{12}^2+J_{23}^2+J_{31}^2\right) -\left( J_{10}^2+J_{20}^2+J_{30}^2\right) $
and $J_{12}J_{30}+J_{23}J_{10}+J_{31}J_{20}$, respectively. $\Omega _{11}$
is an analogue of $J_{12}.$ However we do not know what the analogues of
pure boosts are. There are several candidates that go to $J_{0i}$ in the $%
\lambda \rightarrow 0$ limit and generate a closed algebra with $\Omega ,$
but none of them generates a closed algebra with $P$ and $\Omega $, i.e. $%
\Gamma $ and $\overline{\Gamma }$ are necessary to close the algebra.

\section{Deformed Euclidean subalgebra\label{eu4}}

We denote the subalgebra generated by all $P^{\prime }s$ and all $\Omega
^{\prime }s$ by $E_3^{\lambda ,\hbar }.$ It has the following Casimirs: 
\begin{eqnarray*}
&&Tr_q(P),\qquad \  \\
&&(P,P)_q, \\
&&Tr_q(P\Omega )\stackunder{\lambda \rightarrow 0}{\rightarrow }%
-2P_0+2\lambda \left( 2\overrightarrow{P}\overrightarrow{J}-\hbar \left(
P_0+P_3\right) \right) .
\end{eqnarray*}
The first Casimir is an analogue of the energy and by can be considered to
be constant for an irreducible representation of the algebra. In addition to
the Casimirs one can add $Tr_q(\Omega )$ and $\Omega _{11}$ to form a
complete set of commuting operators for $E_3^{\lambda ,\hbar }$.

\section{Complete set of mutually commuting operators for the full algebra%
\label{set4}}

Below we give a complete set of mutually commuting operators for our
deformed Poincar\'e algebra: 
\begin{eqnarray}
\mathcal{C}_1 &=&(P,P)_q~,\qquad \mathcal{C}_2=(W,W)_q,  \nonumber \\
\mathcal{K}_1 &=&Tr_q(P),  \nonumber \\
\mathcal{K}_2 &=&Tr_q(W)=a\left( \frac \beta {q^3}Tr_q(P\Omega )-\mathcal{K}%
_1\right) ,  \label{pset4} \\
\mathcal{K}_3 &=&Tr_q(\Omega ),\qquad \mathcal{K}_4=\Omega _{11}.  \nonumber
\end{eqnarray}
In the $\lambda \rightarrow 0$ limit $\mathcal{K}_1$ and $\mathcal{K}_2$ go
to $-2P_0$ ($\sim $ energy) and $-2W_0$ ($\sim $ projection of 3-momentum on
3-angular momentum) respectively. The limits and physical sense of the other
operators in the set have been already discussed above.

\section{Eigenstates\label{eigen4}}

In this section we give different procedures for generating the eigenstates
of our set of commuting operators.

\subsection{Eigenstates of $Tr_q(\Omega )$ and $\Omega _{11}$ \label{om4}}

Let us start by considering the eigenstates of $Tr_q(\Omega )$ and $\Omega
_{11}$ constructed in section \ref{su24}. If we have a set of eigenstates $%
\Sigma _j=\{\left| j,m\right\rangle ,\quad m=-j..j,\quad Tr_q(\Omega )\left|
j,m\right\rangle =k_j\left| j,m\right\rangle ,\quad \Omega _{11}\left|
j,m\right\rangle =q^{2m}\left| j,m\right\rangle \}$ we can obtain states
with different angular momenta by applying $Z=\Gamma ,$ $\overline{\Gamma },$
$W,$ $P$ to $\Sigma _j$ as described below. There are four different
procedures associated with the four independent operators in $Z.$

\textbf{Procedure 1: }Act with $Z_{21}$ on the highest weight state $\left|
j,j\right\rangle :$%
\begin{equation}
T_1=Z_{21}\left| j,j\right\rangle .  \label{pr14}
\end{equation}
From the commutational relations (\ref{sigmacr4}) we can compute the
eigenvalues of $Tr_q(\Omega )$ and $\Omega _{11}:$%
\[
Tr_q(\Omega )T_1=k_{j+1}T_1,\quad \Omega _{11}T_1=q^{2\left( j+1\right)
}T_1. 
\]

\textbf{Procedure 2: }Act with $Tr_q(Z)$ on the highest weight state $\left|
j,j\right\rangle :$%
\begin{equation}
T_2=Tr_q(Z)\left| j,j\right\rangle .  \label{pr24}
\end{equation}
Then 
\[
Tr_q(\Omega )T_2=k_jT_2,\quad \Omega _{11}T_2=q^{2j}T_2 
\]

\textbf{Procedure 3: }Define (for $j\geq \frac 12$) 
\begin{equation}
T_3=\left( Z_{11}-Z_{22}\right) \left| j,j\right\rangle +\frac{%
q^2(q^2+1)A_{j,j}}{\left( q^{2j}-q^{-2j}\right) }Z_{21}\left|
j,j-1\right\rangle .  \label{pr34}
\end{equation}
Then 
\[
Tr_q(\Omega )T_3=k_jT_3,\quad \Omega _{11}T_3=q^{2j}T_3. 
\]

\textbf{Procedure 4: }Define (for $j\geq 1$) 
\begin{eqnarray}
T_4 &=&\left( Z_{11}-Z_{22}\right) \left| j,j-1\right\rangle +\frac{%
q^3A_{j,j-1}}{\left( q^{\left( 2j-1\right) }-q^{-\left( 2j-1\right) }\right) 
}Z_{21}\left| j,j-2\right\rangle -  \nonumber \\
&&\frac{q^2B_{j,j-1}}{(q^2-1)q^{2j}}Z_{12}\left| j,j\right\rangle .
\label{pr44}
\end{eqnarray}
Then 
\[
Tr_q(\Omega )T_4=k_{j-1}T_4,\quad \Omega _{11}T_4=q^{2\left( j-1\right)
}T_4. 
\]

Other states with fixed $j$ but different eigenvalues of $\Omega _{11}$ can
be found by acting with $\Omega _{12}$ on the highest weight $T_i$. We shall
only consider below the eigenstates generated by the above procedures from
the rest state, i.e. a state annihilated by $P_{12},$ $P_{21}$ and $\left(
P_{11}-P_{22}\right) .$ We note that acting with $\Gamma $ and $\overline{%
\Gamma }$ can change the energy because $\left[ \Gamma _{ij},Tr_q(P)\right]
\neq 0,$ $\left[ \overline{\Gamma }_{ij},Tr_q(P)\right] \neq 0$. On the
other hand, acting with $P$ and $W$ will not change the energy because $%
\left[ P_{ij},Tr_q(P)\right] =0,$ $\left[ W_{ij},Tr_q(P)\right] =0$. Because
we are interested in the energy spectrum we shall only consider substituting 
$\Gamma $ and $\overline{\Gamma }$ instead of $Z$ in the above procedures.
As we shall see later, not all of the states obtained this way are
eigenvalues of energy and further diagonalization is necessary. For spin 0
states to have 0 eigenvalue of $\left( W,W\right) _q$ one must also choose a
particular value for $\beta $ appearing in (\ref{w4}).

\subsection{Labeling a state \label{label4}}

Let us label an eigenstate of the operators (\ref{pset4}) by 1) $M$ denoting
the ''mass`` (the eigenvalue of $(P,P)_q$ is $M^2$), 2) $s$ denoting the
''spin`` (by this we mean that the rest state eigenvalue of $Tr_q$ $\left(
\Omega \right) $ is $k_s$), 3) the total angular momentum $j$ (the
eigenvalue of $Tr_q$ $\left( \Omega \right) $), 4) ''projection of total
angular momentum on the third axis'' $m$ (the eigenvalue of $\Omega _{11}$
is $q^{2m}$). A state is thus denoted by $\left| M,s,j,m\right\rangle .$ As
we shall see, eigenvalues of other operators from (\ref{pset4}) can be
expressed via these 4 quantum numbers and also depend on the way in which $%
\Gamma $ and $\overline{\Gamma }$ are applied to the rest state.

\subsection{Rest state \label{rest4}}

Let us consider the highest weight rest state $\left| M,s,s,s\right\rangle $
with spin $s$ and mass $M.$ From the fact that $P_{12}$, $P_{21}$ and $%
\left( P_{11}-P_{22}\right) $ should annihilate this state and that the
eigenvalue of $(P,P)_q$ is $M^2$ one can obtain: 
\begin{eqnarray}
P_{11}\left| M,s,s,s\right\rangle &=&-M\left| M,s,s,s\right\rangle ,\quad
P_{12}\left| M,s,s,s\right\rangle =0,\quad  \label{prest4} \\
P_{21}\left| M,s,s,s\right\rangle &=&0,\quad P_{22}\left|
M,s,s,s\right\rangle =-M\left| M,s,s,s\right\rangle .  \nonumber
\end{eqnarray}
From (\ref{prest4}) and 
\begin{equation}
Tr_q\left( \Omega \right) \left| M,s,s,s\right\rangle =k_s\left|
M,s,s,s\right\rangle ,\quad \Omega _{11}\left| M,s,s,s\right\rangle
=q^{2s}\left| M,s,s,s\right\rangle  \label{omegarest4}
\end{equation}
one can calculate eigenvalues of $Tr_q(W)$ and $\left( W,W\right) _q$ . It
turns out that if we want $\left( W,W\right) _q\left| M,0,0,0\right\rangle
=0 $ and $Tr_q(W)\left| M,0,0,0\right\rangle =0$ as in non-deformed case, we
must fix the value of $\beta $ in (\ref{w4}) to be 
\begin{equation}
\beta =q^3.  \label{beta4}
\end{equation}
It is convenient for calculational purposes to rewrite $Tr_q(W)$ in the
form: 
\begin{eqnarray*}
\frac{Tr_q(W)}a &=&\left( P_{11}-P_{22}\right) \left( \Omega _{11}-\frac{%
Tr_q(\Omega )}{q^2+1}\right) +P_{12}\Omega _{21}+ \\
&&q^2P_{21}\Omega _{12}+Tr_q(P)\left( \frac{Tr_q(\Omega )}{q^2+1}-1\right)
\end{eqnarray*}
and to take into account that $\left( W,W\right) _q=-a\left( q^2+1\right)
\left( W,P\right) _q$ for such a choice of $\beta $ (see \cite{sy96a} for
details). Then 
\begin{eqnarray}
Tr_q(W)\left| M,s,s,s\right\rangle &=&-Ma\left( k_s-k_0\right) \left|
M,s,s,s\right\rangle ,  \label{wrest4} \\
(W,W)_q\left| M,s,s,s\right\rangle &=&-Ma\left( q^2+1\right) Tr_q(W)\left|
M,s,s,s\right\rangle  \nonumber \\
&=&-M^2a^2k_0\left( k_s-k_0\right) \left| M,s,s,s\right\rangle .  \nonumber
\end{eqnarray}
It is now obvious that these eigenvalues vanish for zero spin. The
eigenvalues have the limiting values: 
\begin{eqnarray*}
&&\ -Ma\left( k_s-k_0\right) \stackunder{\lambda \rightarrow 0}{\rightarrow }%
-2M\hbar ^2s\left( s+1\right) \cdot \lambda \stackunder{\lambda \rightarrow 0%
}{\rightarrow }0, \\
&&\ -M^2a^2k_0\left( k_s-k_0\right) \stackunder{\lambda \rightarrow 0}{%
\rightarrow }-M^2\hbar ^2s\left( s+1\right) .
\end{eqnarray*}
Note that if a state is obtained from some rest state by applying $\Omega
_{ij}$ it is still a rest state due to the $\Omega -$ $P$ commutational
relations.

\subsection{Procedure 1 \label{Pr14}}

If we apply a monomial $f$ of $l$ order in $\Gamma _{21}$ and $\overline{%
\Gamma }_{21}$ to $\left| M,s,s,s\right\rangle ,$ we get new eigenstates of
our set of commuting operators with the following eigenvalues: 
\begin{eqnarray}
\Omega _{11}f_l\left( \Gamma _{21},\overline{\Gamma }_{21}\right) \left|
M,s,s,s\right\rangle &=&q^{2(l+s)}f_l\left( \Gamma _{21},\overline{\Gamma }%
_{21}\right) \left| M,s,s,s\right\rangle ,  \label{om11p14} \\
Tr_q\left( \Omega \right) f_l\left( \Gamma _{21},\overline{\Gamma }%
_{21}\right) \left| M,s,s,s\right\rangle &=&k_{l+s}f_l\left( \Gamma _{21},%
\overline{\Gamma }_{21}\right) \left| M,s,s,s\right\rangle ,
\label{tromegap14} \\
Tr_q\left( P\right) f_l\left( \Gamma _{21},\overline{\Gamma }_{21}\right)
\left| M,s,s,s\right\rangle &=&-Mk_{\frac l2}f_l\left( \Gamma _{21},%
\overline{\Gamma }_{21}\right) \left| M,s,s,s\right\rangle ,  \label{tpp14}
\\
Tr_q\left( W\right) f_l\left( \Gamma _{21},\overline{\Gamma }_{21}\right)
\left| M,s,s,s\right\rangle &=&-aM\left( k_{s+\frac l2}-k_{\frac l2}\right)
\label{twp14} \\
&&\ \ f_l\left( \Gamma _{21},\overline{\Gamma }_{21}\right) \left|
M,s,s,s\right\rangle ,  \nonumber \\
\left( P_{11}-P_{22}\right) f_l\left( \Gamma _{21},\overline{\Gamma }%
_{21}\right) \left| M,s,s,s\right\rangle &=&M\left( q^l-q^{-l}\right)
\label{p11mp22p14} \\
&&\ \ f_l\left( \Gamma _{21},\overline{\Gamma }_{21}\right) \left|
M,s,s,s\right\rangle .  \nonumber
\end{eqnarray}
Notice from (\ref{om11p14}) and (\ref{tromegap14}) that $l$ and $s$ add like
spin and orbital momentum, therefore we shall interpret $l$ as quantum
number of orbital excitations produced by $f_l\left( \Gamma _{21},\overline{%
\Gamma }_{21}\right) $. Also notice from (\ref{tpp14}) that the energy 
\begin{equation}
E_l=\frac{Mk_{\frac l2}}{k_0}  \label{repenergy}
\end{equation}
depends only on $l$ and is discrete. In the limit $\lambda \rightarrow 0$
the energy $\rightarrow M$. As in non-deformed case $Tr_q\left( W\right) $
(an analogue of $\left( \overrightarrow{J},\overrightarrow{P}\right) $ )
gives $0$ when applied to spin $0$ states for any value of $l.$ For
arbitrary spin $s$ the eigenvalue of $Tr_q\left( W\right) $ goes to 0 when $%
\lambda \rightarrow 0.$ The $\left( P_{11}-P_{22}\right) $ operator (which
corresponds to the projection of the momentum on the third axis) is not, of
course, in the set of commuting operators, only highest weight module is its
eigenstate. In non-deformed case it has $0$ eigenvalue, here it does not.
Note also that the eigenvalues are the same for any internal structure of
the monomial. There seems to be no way to distinguish, for instance, $\Gamma
_{21}\left| M,s,s,s\right\rangle $ from $\overline{\Gamma }_{21}\left|
M,s,s,s\right\rangle $ by eigenvalues of the observables although $\Gamma $
and $\overline{\Gamma }$ commute differently with other elements of the
algebra and themselves.

Since 
\[
\left[ \Omega _{ij},Tr_q(Z)\right] =0, 
\]
states with different eigenvalues of $\Omega _{11}$, obtained by applying
the lowering operator $\Omega _{12}$ to highest weight states, have the same
eigenvalues of $Tr_q(\Omega ),$ $Tr_q(P)$ and $Tr_q(W).$ Of course, they
also have the same eigenvalues of Casimirs $\left( P,P\right) _q$ and$\
\left( W,W\right) _q.$ The eigenvalues of Casimirs can not be changed by
applying any operator to the rest state and characterize the representation
of the whole algebra.

\subsection{Procedures 2 and 3 \label{Pr234}}

This case is more complicated than the previous one because in general
states diagonal in $Tr_q(\Omega )$ and $\Omega _{11}$ would not be
automatically diagonal in other operators in the complete set and further
diagonalization is necessary. The way diagonalization is performed depends
on spin. Let us consider several cases.

\subsubsection{Spin 0 \label{spin04}}

In this case procedure 3 is not applicable. We can only consider procedure 2.

It can be shown that $\left( Tr_q(\Gamma )\right) ^n\left|
M,0,0,0\right\rangle $ for $n$ greater than 1 is not an eigenstate of $%
Tr_q(P).$ Instead the following polynomials of $Tr_q(P)$ applied to $\left|
M,0,0,0\right\rangle $ are eigenstates of energy: 
\begin{eqnarray*}
\pi _0 &=&\left| M,0,0,0\right\rangle ,\quad \\
Tr_q(P)\pi _0 &=&-Mk_0\pi _0, \\
\pi _1 &=&Tr_q(\Gamma )\left| M,0,0,0\right\rangle ,\quad \\
Tr_q(P)\pi _1 &=&-Mk_{\frac 12}\pi _1, \\
\pi _2 &=&\left( \left( Tr_q(\Gamma )\right) ^2-q^2\right) \left|
M,0,0,0\right\rangle ,\quad \\
Tr_q(P)\pi _2 &=&-Mk_1\pi _2, \\
\pi _3 &=&\left( \left( Tr_q(\Gamma )\right) ^3-2q^2Tr_q(\Gamma )\right)
\left| M,0,0,0\right\rangle ,\quad \\
Tr_q(P)\pi _3 &=&-Mk_{\frac 32}\pi _3, \\
\pi _4 &=&\left( \left( Tr_q(\Gamma )\right) ^4-3q^2\left( Tr_q(\Gamma
)\right) ^2+q^4\right) \left| M,0,0,0\right\rangle ,\quad \\
Tr_q(P)\pi _4 &=&-Mk_2\pi _4, \\
\pi _5 &=&\left( \left( Tr_q(\Gamma )\right) ^5-4q^2\left( Tr_q(\Gamma
)\right) ^3+3q^4Tr_q(\Gamma )\right) \left| M,0,0,0\right\rangle ,\quad \\
Tr_q(P)\pi _5 &=&-Mk_{\frac 52}\pi _5, \\
&&\ \ etc.
\end{eqnarray*}
All these results can be compactly described by one recursive relation: 
\begin{eqnarray}
\pi _n-Tr_q(\Gamma )\pi _{n-1} &=&-q^2\pi _{n-2},\quad \pi _0=\left|
M,0,0,0\right\rangle ,\quad \pi _1=Tr_q(\Gamma )\left| M,0,0,0\right\rangle ,
\nonumber \\
\ Tr_q(P)\pi _n &=&-Mk_{\frac n2}\pi _n,\quad Tr_q(W)\pi _n=0,\quad \left(
W,W\right) _q\pi _n=0,\quad  \label{pin4} \\
Tr_q(\Omega )\pi _n &=&k_0\pi _n,\quad \Omega _{11}\pi _n=\pi _n.  \nonumber
\end{eqnarray}
$\pi _n$ is not an eigenstate of $\left( P_{11}-P_{22}\right) .$

If one substitutes $\overline{\Gamma }$ instead of $\Gamma $ everywhere in (%
\ref{pin4}), the eigenvalues are unchanged. However, mixing $\overline{%
\Gamma }$ and $\Gamma $ is much more difficult than in the case of procedure
1 since $\pi _n$ is not a monomial, besides one must take into account that
some quadratic combinations of $\overline{\Gamma }$ and $\Gamma $ should be
converted into 1 or $\Omega .$ We do not discuss this problem here.

\subsubsection{Spin $1/2$ \label{spinhalf4}}

In this case procedures 2 and 3 have to be mixed and it gets even more
complicated. We will consider only the first order in $\Gamma $ and $%
\overline{\Gamma }$ case. The following are eigenstates of all the operators
in the complete set of commuting operators: 
\begin{eqnarray}
S_1 &=&Tr_q\left( \Gamma \right) \left| M,\frac 12,\frac 12,\frac
12\right\rangle +  \nonumber \\
&&q^3\left( \Gamma _{21}\left| M,\frac 12,\frac 12,-\frac 12\right\rangle
-\frac 1q\Gamma _{22}\left| M,\frac 12,\frac 12,-\frac 12\right\rangle
\right) ,  \nonumber \\
S_2 &=&Tr_q\left( \Gamma \right) \left| M,\frac 12,\frac 12,\frac
12\right\rangle -  \nonumber \\
&&q^5\left( \Gamma _{21}\left| M,\frac 12,\frac 12,-\frac 12\right\rangle
-\frac 1q\Gamma _{22}\left| M,\frac 12,\frac 12,-\frac 12\right\rangle
\right) ,  \nonumber \\
S_3 &=&Tr_q\left( \overline{\Gamma }\right) \left| M,\frac 12,\frac 12,\frac
12\right\rangle +  \label{Si4} \\
&&q\left( \overline{\Gamma }_{21}\left| M,\frac 12,\frac 12,-\frac
12\right\rangle -\frac 1q\overline{\Gamma }_{22}\left| M,\frac 12,\frac
12,-\frac 12\right\rangle \right) ,  \nonumber \\
S_4 &=&Tr_q\left( \overline{\Gamma }\right) \left| M,\frac 12,\frac 12,\frac
12\right\rangle -  \nonumber \\
&&\frac 1q\left( \overline{\Gamma }_{21}\left| M,\frac 12,\frac 12,-\frac
12\right\rangle -\frac 1q\overline{\Gamma }_{22}\left| M,\frac 12,\frac
12,-\frac 12\right\rangle \right)  \nonumber
\end{eqnarray}
with eigenvalues 
\begin{eqnarray}
\Omega _{11}S_i &=&qS_i,\quad Tr_q\left( \Omega \right) S_i=k_{\frac
12}S_i,\quad Tr_q\left( P\right) S_i=-mk_{\frac 12}S_i;  \nonumber \\
Tr_q\left( W\right) S_{1,3} &=&aM\left( k_{\frac 12}-k_0\right) S_{1,3},\quad
\label{Ei4} \\
Tr_q\left( W\right) S_{2,4} &=&-aM\left( k_1-k_{\frac 12}\right) S_{1,4}, 
\nonumber \\
\left( P_{11}-P_{22}\right) S_{1,3} &=&M\left( q-\frac 1q\right)
S_{1,3},\quad  \nonumber \\
&&S_{2,4}\quad is\quad not\quad eigenstate\quad of\quad \left(
P_{11}-P_{22}\right) .  \nonumber
\end{eqnarray}

We do not know the general formula for higher orders here. Notice again that
one can not distinguish here states generated by $\Gamma $ and $\overline{%
\Gamma }:$ $S_{1,3}$ have the same eigenvalues for any observables, so do $%
S_{2,4};$ but $S_3$ $\left( S_4\right) $ is not obtained by just
substituting $\overline{\Gamma }$ instead of $\Gamma $ into $S_1$ $\left(
S_2\right) :$ the mixing coefficients are different.

\subsection{Procedure 4 \label{Pr44}}

We only consider here how to use procedure 4 to shift from the rest state
with spin 1 to the excited state with angular momentum 0. 
\begin{eqnarray}
S_5 &=&-\left( \Gamma _{11}-\Gamma _{22}\right) \left| M,1,1,0\right\rangle +
\nonumber \\
&&\ \ \frac{\sqrt{q^2+1}}{q^2}\left( \Gamma _{12}\left| M,1,1,1\right\rangle
-q^3\Gamma _{21}\left| M,1,1,-1\right\rangle \right) ,  \nonumber \\
\Omega _{11}S_5 &=&S_5,\quad Tr_q\left( \Omega \right) S_5=k_0S_5,\quad
\label{P44} \\
Tr_q\left( P\right) S_5 &=&-Mk_{\frac 12}S_5,\quad Tr_q\left( W\right) S_5=0.
\nonumber
\end{eqnarray}
Again, we can substitute $\overline{\Gamma }$ instead of $\Gamma $ above.

\section{Conclusion\label{concl4}}

Eigenvalues for all the states found above in the limit $\lambda \rightarrow
0$ go to those of the rest state of the non-deformed Poincar\'e algebra.
Therefore the deformation splits the rest state into an infinite number of
states with a discrete energy spectrum which is bounded from below by $M$
but is unbounded from above. Within this interpretation, the observed
curious strong degeneracy of states (when states constructed by applying
algebraically independent operators to the rest state are nevertheless
indistinguishable by eigenvalues of all the observables in the algebra) is
not that surprising since in the limit they all have the same eigenvalues.

States generated by procedures 1 and 4 have different energies and angular
momenta and hence can be thought of as ''quantum rotations'', while states
generated by procedures 2 and 3 have different energies but same angular
momentum and therefore can be interpreted as ''quantum oscillations'' which
a particle acquires upon the deformation.

It is remarkable that both ''rotations'' and ''oscillations'' have exactly
the same expression $\left( \text{\ref{repenergy}}\right) $ for the energy
spectrum (though different quantum numbers, $l$ and $n,$ are substituted
into this expression).

It would be interesting to see what happens to a moving particle and how $E=%
\sqrt{P^2+M^2}$ is modified. To answer this question, one must either figure
out how to apply deformed Lorentz transformations to these ''rest'' states
(and also find an analogue of boosts) or start with a different set of
commuting operators: in \cite{sy96a} we found an alternative set of
commuting operators in which $Tr_q\left( \Omega \right) $ is exchanged for $%
P_3=P_{11}-P_{22}.$ Such a set could probably be also used to construct
representations of a massless particle\footnote{%
For a massless particle we must also start from a different ground state
since there is no rest state in this case: 
\[
P_{12}\left| 0\right\rangle =P_{21}\left| 0\right\rangle =P_{11}\left|
0\right\rangle =0, 
\]
or 
\[
P_{12}\left| 0\right\rangle =P_{21}\left| 0\right\rangle =P_{22}\left|
0\right\rangle =0. 
\]
}.

Of course, one might also try to find general formulas for the states
presented in this paper: we found some but could not find others. Also we
did not consider here what happens when one mixes ''rotations'' with
''oscillations'' or even different kinds of ''rotations'' (obtained using
procedures 1 and 4). However it was very computationally difficult and
required enormous computer resources to find the results we obtained despite
the fact that the results themselves are surprisingly simple and nice.
Perhaps there might exist such a point of view from which all these results
are obvious and do not require such difficult calculations.

Another problem which is left open is the question of adding representations
for two or more relativistic particles. The answer should be nontrivial
since $\Gamma $, $\overline{\Gamma }$ and $P$ does not appear to generate a
Hopf algebra.

\section{Acknowledgment}

I am grateful to Dr. A. Stern for useful discussions and his help. This work
was in part supported by Graduate Council Research Fellowship from the
University of Alabama and by Department of Energy, USA, under contract
number DE-FG0584ER40141.

\bibliographystyle{h-physre}
\bibliography{qg1}

\section{Appendix: commutational relations (\ref{pggb4}) in components\label
{app4}}

\subsection{$P-P$}

\begin{eqnarray*}
P_{12}~P_{11} &=&P_{11}~P_{12}+(1-q^2)~P_{12}~P_{22},  \label{p11p124} \\
P_{11}~P_{21} &=&P_{21}~P_{11}+(1-q^2)~P_{22}~P_{21}, \\
P_{11}~P_{22} &=&P_{22~}P_{11}, \\
P_{12}~P_{21} &=&P_{21}~P_{12}+(q^2-1)~P_{22}~(P_{11}-P_{22}), \\
P_{12}~P_{22} &=&q^2P_{22}~P_{12}, \\
P_{22}~P_{21} &=&q^2P_{21}~P_{22}.
\end{eqnarray*}

\subsection{$\Gamma -\Gamma $}

\begin{eqnarray*}
\Gamma _{12}\Gamma _{11} &=&q^2\Gamma _{11}\Gamma _{12},  \label{g11g124} \\
\Gamma _{21}\Gamma _{11} &=&\frac 1{q^2}\Gamma _{11}\Gamma _{21}, \\
\Gamma _{22}\Gamma _{11} &=&\Gamma _{11}\Gamma _{22}, \\
\Gamma _{21}\Gamma _{12} &=&\Gamma _{12}\Gamma _{21}+\frac{1-q^2}{q^2}\Gamma
_{11}(\Gamma _{22}-\Gamma _{11}), \\
\Gamma _{22}\Gamma _{12} &=&\Gamma _{12}\Gamma _{22}-\frac{1-q^2}{q^2}\Gamma
_{11}\Gamma _{12}, \\
\Gamma _{22}\Gamma _{21} &=&\Gamma _{21}\Gamma _{22}+\frac{1-q^2}{q^2}\Gamma
_{21}\Gamma _{11}.
\end{eqnarray*}

\subsection{$\overline{\Gamma }-\overline{\Gamma }$}

\begin{eqnarray*}
\overline{\Gamma }_{12}\overline{\Gamma }_{11} &=&\overline{\Gamma }_{11}%
\overline{\Gamma }_{12}+(1-q^2)\overline{\Gamma }_{12}\overline{\Gamma }%
_{22},  \label{gb11gb124} \\
\overline{\Gamma }_{21}\overline{\Gamma }_{11} &=&\overline{\Gamma }_{11}%
\overline{\Gamma }_{21}-(1-q^2)\overline{\Gamma }_{22}\overline{\Gamma }%
_{21}, \\
\overline{\Gamma }_{22}\overline{\Gamma }_{11} &=&\overline{\Gamma }_{11}%
\overline{\Gamma }_{22}, \\
\overline{\Gamma }_{21}\overline{\Gamma }_{12} &=&\overline{\Gamma }_{12}%
\overline{\Gamma }_{21}-(1-q^2)\overline{\Gamma }_{22}(\overline{\Gamma }%
_{22}-\overline{\Gamma }_{11}), \\
\overline{\Gamma }_{22}\overline{\Gamma }_{12} &=&\frac 1{q^2}\overline{%
\Gamma }_{12}\overline{\Gamma }_{22}, \\
\overline{\Gamma }_{22}\overline{\Gamma }_{21} &=&q^2\overline{\Gamma }_{21}%
\overline{\Gamma }_{22}.
\end{eqnarray*}

\subsection{$\Gamma -\overline{\Gamma }$}

\begin{eqnarray*}
\overline{\Gamma }_{11}\Gamma _{11} &=&\Gamma _{11}\overline{\Gamma }%
_{11}+(1-q^2)\Gamma _{12}\overline{\Gamma }_{21},  \label{g11gb114} \\
\overline{\Gamma }_{12}\Gamma _{11} &=&\Gamma _{11}\overline{\Gamma }%
_{12}+(1-q^2)\Gamma _{12}\overline{\Gamma }_{22}-(1-q^2)\overline{\Gamma }%
_{11}\Gamma _{12}, \\
\overline{\Gamma }_{11}\Gamma _{12} &=&\Gamma _{12}\overline{\Gamma }_{11},
\\
\overline{\Gamma }_{12}\Gamma _{12} &=&q^2\Gamma _{12}\overline{\Gamma }%
_{12}, \\
\overline{\Gamma }_{21}\Gamma _{11} &=&\Gamma _{11}\overline{\Gamma }_{21},
\\
\overline{\Gamma }_{22}\Gamma _{11} &=&\Gamma _{11}\overline{\Gamma }%
_{22}-(1-q^2)\overline{\Gamma }_{21}\Gamma _{12}, \\
\overline{\Gamma }_{21}\Gamma _{12} &=&\frac 1{q^2}\Gamma _{12}\overline{%
\Gamma }_{21}, \\
\overline{\Gamma }_{22}\Gamma _{12} &=&\Gamma _{12}\overline{\Gamma }_{22},
\\
\overline{\Gamma }_{11}\Gamma _{21} &=&\Gamma _{21}\overline{\Gamma }%
_{11}+(1-q^2)(\Gamma _{22}-\Gamma _{11})\overline{\Gamma }_{21}, \\
\overline{\Gamma }_{12}\Gamma _{21} &=&\frac 1{q^2}\Gamma _{21}\overline{%
\Gamma }_{12}+\frac{(1-q^2)}{q^2}(\Gamma _{22}-\Gamma _{11})\overline{\Gamma 
}_{22}- \\
&&\frac{(1-q^2)}{q^2}\overline{\Gamma }_{11}(\Gamma _{22}-\Gamma _{11}), \\
\overline{\Gamma }_{11}\Gamma _{22} &=&\Gamma _{22}\overline{\Gamma }_{11}-%
\frac{(1-q^2)}{q^2}\Gamma _{12}\overline{\Gamma }_{21}, \\
\overline{\Gamma }_{12}\Gamma _{22} &=&\Gamma _{22}\overline{\Gamma }_{12}-%
\frac{(1-q^2)}{q^2}\Gamma _{12}\overline{\Gamma }_{22}+\ \frac{(1-q^2)}{q^2}%
\overline{\Gamma }_{11}\Gamma _{12}, \\
\overline{\Gamma }_{12}\Gamma _{22} &=&\Gamma _{22}\overline{\Gamma }_{12}-%
\frac{(1-q^2)}{q^2}\Gamma _{12}\overline{\Gamma }_{22}+\ \frac{(1-q^2)}{q^2}%
\overline{\Gamma }_{11}\Gamma _{12}, \\
\overline{\Gamma }_{21}\Gamma _{21} &=&q^2\Gamma _{21}\overline{\Gamma }%
_{21}, \\
\overline{\Gamma }_{22}\Gamma _{21} &=&\Gamma _{21}\overline{\Gamma }_{22}-\ 
\frac{(1-q^2)}{q^2}\overline{\Gamma }_{21}(\Gamma _{22}-\Gamma _{11}), \\
\overline{\Gamma }_{21}\Gamma _{22} &=&\Gamma _{22}\overline{\Gamma }_{21},
\\
\overline{\Gamma }_{22}\Gamma _{22} &=&\Gamma _{22}\overline{\Gamma }_{22}+%
\frac{(1-q^2)}{q^2}\overline{\Gamma }_{21}\Gamma _{12}.
\end{eqnarray*}

\subsection{$P-\Gamma $}

\begin{eqnarray*}
\Gamma _{11}P_{11} &=&qP_{11}\Gamma _{11}-(1-q^2)\Gamma _{12}P_{21},
\label{p11g114} \\
\Gamma _{12}P_{11} &=&\frac 1q\ P_{11}\Gamma _{12}-\frac{(1-q^2)}{q^2}\Gamma
_{11}P_{12}-\frac{(1-q^2)^2}{q^2}\Gamma _{12}P_{22}, \\
\Gamma _{11}P_{12} &=&\frac 1qP_{12}\Gamma _{11}-(1-q^2)\Gamma _{12}P_{22},
\\
\Gamma _{12}P_{12} &=&\frac 1qP_{12}\Gamma _{12}, \\
\Gamma _{21}P_{11} &=&qP_{11}\Gamma _{21}-(1-q^2)\Gamma
_{22}P_{21}+q(1-q^2)P_{21}\Gamma _{11}, \\
\Gamma _{22}P_{11} &=&\frac 1qP_{11}\Gamma _{22}-\frac{(1-q^2)}{q^2}\Gamma
_{21}P_{12}-\frac{(1-q^2)^2}{q^2}\Gamma _{22}P_{22}+ \\
&&\frac{(1-q^2)}qP_{21}\Gamma _{12}, \\
\Gamma _{21}P_{12} &=&\ qP_{12}\Gamma _{21}-(1-q^2)\Gamma _{22}P_{22}+\frac{%
(1-q^2)}q(P_{22}-P_{11})\Gamma _{11}, \\
\Gamma _{22}P_{12} &=&qP_{12}\Gamma _{22}+\frac{(1-q^2)}q(P_{22}-P_{11})%
\Gamma _{12}, \\
\Gamma _{11}P_{21} &=&qP_{21}\Gamma _{11}, \\
\Gamma _{12}P_{21} &=&qP_{21}\Gamma _{12}-(1-q^2)\Gamma _{11}P_{22}, \\
\Gamma _{11}P_{22} &=&\frac 1qP_{22}\Gamma _{11}, \\
\Gamma _{12}P_{22} &=&qP_{22}\Gamma _{12}, \\
\Gamma _{21}P_{21} &=&\frac 1qP_{21}\Gamma _{21}, \\
\Gamma _{22}P_{21} &=&\frac 1qP_{21}\Gamma _{22}-(1-q^2)\Gamma _{21}P_{22},
\\
\Gamma _{21}P_{22} &=&\frac 1qP_{22}\Gamma _{21}-\frac{(1-q^2)}{q^3}%
P_{21}\Gamma _{11}, \\
\Gamma _{22}P_{22} &=&qP_{22}\Gamma _{22}-\frac{(1-q^2)}qP_{21}\Gamma _{12}.
\end{eqnarray*}

\subsection{$P-\overline{\Gamma }$}

\begin{eqnarray*}
\overline{\Gamma }_{11}P_{11} &=&qP_{11}\overline{\Gamma }_{11}-(1-q^2)%
\overline{\Gamma }_{12}P_{21}+q(1-q^2)P_{12}\overline{\Gamma }_{21},
\label{p11gb114} \\
\overline{\Gamma }_{12}P_{11} &=&\frac 1qP_{11}\overline{\Gamma }_{12}+\frac{%
(1-q^2)}qP_{12}(\overline{\Gamma }_{22}-\overline{\Gamma }_{11}), \\
\overline{\Gamma }_{11}P_{12} &=&qP_{12}\overline{\Gamma }_{11}-(1-q^2)%
\overline{\Gamma }_{12}P_{22}, \\
\overline{\Gamma }_{12}P_{12} &=&qP_{12}\overline{\Gamma }_{12}, \\
\overline{\Gamma }_{21}P_{11} &=&qP_{11}\overline{\Gamma }_{21}-(1-q^2)%
\overline{\Gamma }_{22}P_{21}, \\
\overline{\Gamma }_{22}P_{11} &=&\frac 1qP_{11}\overline{\Gamma }_{22}-\frac{%
(1-q^2)}{q^3}P_{12}\overline{\Gamma }_{21}, \\
\overline{\Gamma }_{21}P_{12} &=&\frac 1qP_{12}\overline{\Gamma }%
_{21}-(1-q^2)\overline{\Gamma }_{22}P_{22}, \\
\overline{\Gamma }_{22}P_{12} &=&\frac 1qP_{12}\overline{\Gamma }_{22}, \\
\overline{\Gamma }_{11}P_{21} &=&\frac 1q\ P_{21}\overline{\Gamma }_{11}-%
\frac{(1-q^2)}{q^2}\overline{\Gamma }_{21}P_{11}-\frac{(1-q^2)^2}{q^2}%
\overline{\Gamma }_{22}P_{21}+ \\
&&\frac{(1-q^2)}qP_{22}\overline{\Gamma }_{21}, \\
\overline{\Gamma }_{12}P_{21} &=&\frac 1qP_{21}\overline{\Gamma }%
_{12}-(1-q^2)\overline{\Gamma }_{22}P_{11}+\frac{(1-q^2)}qP_{22}(\overline{%
\Gamma }_{22}-\overline{\Gamma }_{11}), \\
\overline{\Gamma }_{11}P_{22} &=&\frac 1q\ P_{22}\overline{\Gamma }_{11}-%
\frac{(1-q^2)}{q^2}\overline{\Gamma }_{21}P_{12}-\frac{(1-q^2)^2}{q^2}%
\overline{\Gamma }_{22}P_{22}, \\
\overline{\Gamma }_{12}P_{22} &=&qP_{22}\overline{\Gamma }_{12}-(1-q^2)%
\overline{\Gamma }_{22}P_{12}, \\
\overline{\Gamma }_{21}P_{21} &=&qP_{21}\overline{\Gamma }_{21}, \\
\overline{\Gamma }_{22}P_{21} &=&qP_{21}\overline{\Gamma }_{22}-\frac{(1-q^2)%
}qP_{22}\overline{\Gamma }_{21}, \\
\overline{\Gamma }_{21}P_{22} &=&\frac 1qP_{22}\overline{\Gamma }_{21}, \\
\overline{\Gamma }_{22}P_{22} &=&qP_{22}\overline{\Gamma }_{22}.
\end{eqnarray*}

\end{document}